\documentclass[12pt]{amsart}
\usepackage{amsfonts,amsmath,amstext,amsthm,amssymb}
\usepackage{hyperref}
\numberwithin{equation}{section}
\newtheorem{theorem}{Theorem}
\newtheorem{lemma}{Lemma}
\newtheorem{property}{Property}

\def\rank{{\rm rank\,}}
\def\seg{{\rm Int\,}}
\def\dist{{\rm dist\,}}
\def\supp{{\rm supp\,}}
\def\ZR{\ensuremath{\mathbb R}}

\def\ZI{\ensuremath{\mathbb I}}

\def\md#1#2\emd{\ifx0#1
\begin{equation*} #2 \end{equation*}\fi  
\ifx1#1\begin{equation}#2\end{equation}\fi   
\ifx2#1\begin{align*}#2\end{align*}\fi   
\ifx3#1\begin{align}#2\end{align}\fi    
\ifx4#1\begin{gather*}#2\end{gather*}\fi  
\ifx5#1\begin{gather}#2\end{gather}\fi   
\ifx6#1\begin{multline*}#2\end{multline*}\fi  
\ifx7#1\begin{multline}#2\end{multline}\fi  
}
\newcommand {\e }[1]{(\ref{#1})}
\newcommand {\lem }[1]{Lemma \ref{#1}}

\begin{document}
\title{On directional maximal operators\\ associated with
 generalized lacunary sets}

\author
{Grigor A. Karagulyan}

\address{Institute of Mathematics \\
Armenian National Academy of Sciences\\
Marshal Baghramian ave. 24b,\\
Yerevan, 375019, ARMENIA\\
}

\email{karagul@instmath.sci.am}

\date{}
\maketitle

\begin{section}{Introduction}\everypar{} \parskip=11pt

For any number $s\in [0,1]$ we denote $e_s=(\cos 2\pi s ,\sin 2\pi
s) $. We define the following operators on $\ZR^2$ associated with
this direction $s$ by
\md4
M_s^0f(x)=\frac{1}{2}\int_{-1}^1|f(x+te_s)|dt,\\
M_s^1f(x)=\sup_{\delta
>0}\frac{1}{2\delta }\int_{-\delta}^\delta |f(x+te_s)|dt,\\
M_s^2f(x)=\sup_{R\in \mathcal{R}_s, x\in R } \frac{1}{|R|}\int_R
|f(x)|dx
\emd
where $\mathcal{R}_s$ is the set of all rectangles in $\ZR^2$
having slope $s$. For any set of slopes $\Omega $ we denote
\md1\label{MOm}
M_\Omega^j f(x)=\sup_{s\in \Omega}M_s^jf(x),\quad j=0,1,2.
\emd
 Obviously we have
\md1\label{MEq}
M_\Omega^0f(x)\le M_\Omega^1f(x)\le M_\Omega^2f(x)\le
 M_\Omega^1M_{\Omega^\perp }^1f(x).
\emd
where $\Omega^\perp $ is the set of vectors orthogonal some vector
from $\Omega$. The problem of boundedness of these operators for
different $\Omega $ has a long history. The case of lacunary set
$\Omega =\{s_k\}$, $s_{k+1}\le \lambda s_k(\lambda<1)$ is
investigated in the papers J.~Str{\"o}mberg \cite{Str1},
A.~C\'{o}rdoba, R.Fefferman \cite{CoFe}, A.~Nagel, E.M.~Stein and
S.~Wainger \cite{NSW}. A final result is obtained in \cite{NSW} in
1979, where it is proved the boundedness of these operators in
spaces $L^p$, $1<p<\infty $.

The next rank of results concerns the case when $\Omega $ is
finite set of cardinality $N$. The earliest results of related
operators was carried out by A.~C\'{o}rdoba \cite{Cor}. He
obtained a bound $C\sqrt{\log N} $ on $L^2$ norm of maximal
operator $M_\Omega^0f(x)$ in the case of uniformly distributed
$\Omega $. A bound $C\log N$ for the operator $M_\Omega^2f(x)$
for the similar set $\Omega $ is obtained by J.~Str{\"o}mberg in
\cite{Str2}. An interesting problem was if C\'{o}rdoba's or
Str{\"o}mberg's results is extended to the case of $N$ distinct
directions.  A partial result was treated by Barrionuevo
\cite{Barr} , \cite{Barr1}.  And the definitive estimates
\md5
\|M_\Omega^0f(x)\|_{L^2}\le C\sqrt{\log N}\|f\|_{L^2},\label{M0}\\
\|M_\Omega^2f(x)\|_{L^2}\le C\log N\|f\|_{L^2},\label{M2}
\emd
was obtained by N.~Katz in \cite{Katz2}, \cite{Katz1} in 1999.

We are interested in extensions of lacunary sets of directions, to
collections we call $\mu $--lacunary, for an integer $\mu$. We say
interval $J=[a,b]$ is adjacent for the set $A\subset \ZR $ iff
$(a,b)\cap A=\varnothing $ and for any other interval
$(c,d)\supset (a,b)$ we have $(c,d)\cap A\neq \varnothing $.  We
define $\mu $--lacunary set $\Omega $ by induction. Say that the
sequence $\Omega_1=\{v_i \}$ (finite or infinite) is lacunary (or
$1$--lacunary) iff there is a number $v_\infty$ so that
\md1\label{deflac}
|v_{i+1}-v_\infty |<\lambda |v_i-v_\infty |,\quad  0<\lambda <1.
\emd
Every $k+1$--lacunary set can be obtained from some $k$--lacunary
$\Omega_k$ adding some points to $\Omega_k$ as follow. In each
interval adjacent for  $\Omega_k$ we can add a lacunary sequence
(finite or infinite). So if $\Omega$ is some $\mu$--lacunary set
we can fix a sequence of sets
$\Omega_1\subset\Omega_2\subset\cdots
\subset\Omega_{\mu-1}\subset\Omega_\mu=\Omega  $ such that each
$\Omega_k$ is $k$--lacunary.

The boundedness of maximal operator $M_\Omega^2f(x) $ in $L^p,\,
(p>1)$ for $\mu$-lacunary $\Omega $ is proved in P.Sjogren and
P.Sjolin \cite{SS}. We are interested in growth of the norms of
$M_\Omega^1 f(x)$ in $L^2$ for a $\mu$--lacunary $\Omega $, as
$\mu$ tends to infinity. Recently a sharp estimate for the maximal
function $M_\Omega^2 f(x)$ is established in  A.~Alfonesca,
F.~Soria and A.~Vargas \cite{ASV2}, and G.A.Karagulyan, M.Lacey
\cite{KaLa}. That is
\md1\label{M2lac}
\|M_\Omega^2 f(x)\|_{L^2}\lesssim{} \mu \|f\|_{L^2}.
\emd
for any $\mu$--lacunary set $\Omega $. Actually, in A.~Alfonesca,
F.~Soria and A.~Vargas \cite{ASV1}, \cite{ASV2} and A.~Alfonesca
\cite{Alf} it has been obtained a general result, an interesting
orthogonality principle for the maximal function $M_\Omega^2$.  In
this paper we shall prove the better estimate for the maximal
operator $M_\Omega^1 f(x)$.


\begin{theorem}  \label{Tlac}
Let $\Omega$ be any $\mu$--lacunary set. Then we have
\md1\label{M1lac}
\|M_\Omega^1 f(x)\|_{L^2}\lesssim{}\sqrt {\mu}  \|f\|_{L^2}.
\emd
\end{theorem}
According to a simple \lem{L2} below any set $\Omega $ of
cardinality $N$ is $\big([\log N]+2\big)$-lacunary. So the
inequality \e{M2lac} implies \e{M2} and from \e{M1lac} we get a
new estimate
\md1\label{M1}
\|M_\Omega^1f(x)\|_{L^2}\le C\sqrt{\log N}\|f\|_{L^2},
\emd
for any set $\Omega $ of cardinality $N$. Taking account of
\e{MEq}, \e{M1} implies both estimates \e{M0}, \e{M2}. It is known
the sharpness of the orders of constants in \e{M0} and \e{M2}. So
the same we have for \e{M1lac} and \e{M1}.

\end{section}
\begin{section}{Lacunary sets}
Let us go back to the definition of $\mu$--lacunary  set $\Omega
$. In each step of construction we get $k$-lacunary set
$\Omega_k$. We call all adjacent intervals of $\Omega_k$ by
$k$-rank intervals of  $\Omega $. We shall note it by $\rank
(J)=k$. Observe that all intervals of the same rank are mutually
disjoint and if $\rank I>\rank J$ then $I\subset J$ or $I\cap
J=\varnothing $. Denote the set of all rank intervals of $\Omega $
by $\seg\Omega$. We call poles of the $\Omega $ all the points
$v_\infty $ in definition of $\mu $-lacunary set. We connect with
each interval $J\in \seg\Omega $ with $\rank J\le \mu-1$ a pole
$p_J\in J$ which is the first one appeared in the process of
construction.

Now we are going define a complete $\mu $-lacunary set. We
consider a monotone sequence $v=\{v_k:k=1,2,\cdots \}$ satisfying
\md0
\frac{1}{4}|v_k-v_\infty |\le |v_{k+1}-v_\infty |<
\frac{1}{2}|v_k-v_\infty |
\emd
for some $v_\infty $. We say $v=\{v_k\}$ is compete one-side
lacunary in interval $J=(a,b)\supset v$ if $v_1-v_\infty\ge
\frac{1}{2}(b-v_\infty)$ in the case of decreasing $v_k$ and
$v_\infty-v_1\ge \frac{1}{2}(v_\infty-a)$ otherwise. We say a set
is complete both-side lacunary in the interval $J$ if it is union
of two compete one-side lacunary sequences $v$ and $v'$ in $J$
with the same pole $v_\infty $ such that $v$ is increasing and
$v'$ is decreasing. We say a set is just complete lacunary iff it
is complete lacunary in some interval. The following properties
are easy to check.

\begin{property}
If $v$ is a lacunary set with pole $v_\infty$ and gap $\lambda $
then for any $J=[a,b]$ the set $v\cap J$ also is lacunary with gap
$\lambda $ which pole is in $J$.
\end{property}
\begin{property}
Let $v$ be a lacunary set with pole $v_\infty$ and gap $\lambda
<1/2$. If $v\subset J=(a,b)$ and $v_\infty\in J$, then there
exists a complete lacunary set $v^*$ in $J$ such that
$v\subset\overline v$.
\end{property}
\begin{lemma}\label{L1}
Any $\mu $-lacunary set with gap $\lambda $ is a part of some
complete $n(\lambda ) \mu $-lacunary set, where $n(\lambda )=1$ if
$\lambda \le 1/2$ and $n(\lambda )=\log_2^{-1}(1/\lambda )$ if
$1/2<\lambda < 1$.
\end{lemma}
Suppose first $\lambda \le 1/2$. By the definition of $\mu
$-lacunary set $\Omega $ we start the construction of $\Omega $
with an original lacunary set $\Omega_1$. First of all we do any
completion $\Omega_1^*$ of $\Omega_1$. Suppose we already obtained
a complete $k$--lacunary $\Omega_k^*$ containing $\Omega_k$.
Consider all lacunary sets $v=\{v_k\}$ were added to $\Omega_k$ to
get $\Omega_{k+1}$. We note that at most one of such $v=\{v_k\}$
can intersects a fixed adjacent interval $J$ of $\Omega_k^*$. To
get $\Omega_{k+1}^*$ we add in each $J$ a completion of $J\cap v$.
At the end of process we shall get $\mu $-lacunary set
$\Omega^*=\Omega^*_\mu$ containing $\Omega $.

\begin{lemma}\label{L2}
Any set of cardinality $N$ is $\big([\log_2N]+2\big)$--lacunary.
\end{lemma}
\begin{proof}
Let $2^m\le N<2^{m+1}$. We assume $N=2^{m+1}-1$ and $\#\Omega =N$.
Suppose $\Omega =\{x_k:\, k=0,1,\cdots , 2^{m+1}\}$ and
$x_0<x_1<\cdots <x_{2^{m+1}}$. The original lacunary set has to be
taken the set $\Omega_1=\{x_0,x_{2^{m+1}}\}$. In its single
adjacent interval $[x_0,x_{2^{m+1}}]$ we take a lacunary set
consists of the only point $x_{2^m}$. Thus we get
$\Omega_2=\{x_0,x_{2^m},x_{2^{m+1}}\}$. In each adjacent intervals
$[x_0,x_{2^m}]$ and $[x_{2^m},x_{2^{m+1}}]$ we take two lacunary
sets consist of the points $x_{2^{m-1}}$ and $x_{2^m+2^{m-1}}$
correspondingly. It will be the third step of construction. After
the $(m+2)$th step all points of $\Omega $ will be chosen and we
will get an $(m+2)$--lacunary set.
\end{proof}

\end{section}

\begin{section}{Notations}  \everypar{} \parskip=11pt

Remind the Fejer kernel
\md0
K_r(x)=\int_{-r}^r\left(1-\frac{|t|}{r}\right)e^{-itx}dt=
\frac{}{}\frac{4\sin^2\frac{rx}{2}}{rx^2}.
\emd
Vallee-Poussin kernel is defined by
\md0
 V_r(x)=2K_{2r}(x)-K_{r}(x).
\emd
For the Fourier transform of this function we have
\md1\label{lin}
\widehat V_r(\xi)={} \left\{
\begin{array}{rcl}
1 &\hbox{ if }& |\xi|\in [0,r],\\ 0& \hbox{ if }&  |\xi|>2r,\\
\hbox {linear } &\hbox{ on }&\pm [r,2r].
\end{array}
\right.
\emd
From a property of Fejer kernel we have
\md0
| V_r(x)|\le C \max\left\{\frac{1}{rx^2},r\right\}
\emd
Thus numbers $\gamma_k>0$ we get
\md1\label{int}
| V_r(x)|\le C\sum_{k>\log 1/r} \gamma_k\ZI_{(-2^k,2^k)}(x)=
 \zeta_r(x)\in L^1(\ZR ).
\emd
Choose a Schwartz function $\phi $ with
\md1
\phi (x) \ge 0, \quad\phi (x) \ge 1 \hbox { as } x\in [0,1],\quad
\supp\widehat\phi \subset[-1,1]\label{psi}.
\emd
It is clear, that
\md1\label{xi}
\xi (x)=\max\{|\phi (x)|,|x\phi (x)|\}\in L^1(\ZR ).
\emd
Denote $u_\alpha =(1,\alpha )$ consider the directional maximal
operator
\md0
M_\alpha f(x)=\sup_{\delta
>0}\frac{1}{2\delta }\int_{-\delta}^\delta |f(x+tu_\alpha )|dt
\emd
If $\alpha \in[0,1]$ and $s=\frac{\arctan \alpha}{2\pi } $ this
operator is equivalent to the operator $M_s^1f(x)$ defined in the
beginning of the article. Therefore to prove the theorem it is
enough to prove
\md1
\|\sup_{\alpha \in \Omega }M_\alpha f(x)\|_{L^2}\le C\sqrt {\mu}
\|f\|_{L^2}
\emd
for any $\mu $-lacunary set $\Omega\subset [0,1]$. Define
operators
\md1\label{Gam}
\Gamma_{\alpha ,r,h} f(x)=\big(
V_r(x_2-x_1\alpha)\phi_h(x_1)\big)* f(x),\quad x=(x_1,x_2)\in
\ZR^2,
\emd
where
\md0
\phi_h(x)=\frac{1}{h}\phi \left(\frac{x}{h}\right).
\emd
According to the Lebesgue's theorem on differentiation of
integrals we have
\md0
\lim_{h\to 0}\Gamma_{\alpha ,r,h} f(x)=\int_\ZR
V_r(t)f(x+tu_\alpha )dt\hbox { a.e. }.
\emd
On the other hands
\md0
\frac{1}{2r}\ZI_{(-r,r)}(x)\le
K_r(x)=\frac{1}{2}V_{r/2}(x)+\frac{1}{4}V_{r/4}(x)+\cdots
\emd
From this it follows that
\md0
M_\alpha f(x)=\sup_{\delta
>0}\frac{1}{2r}\int_{-r}^r|f(x+tu_\alpha )|dt\le
\sup_{r>0}\left|\int_\ZR V_r(t)|f(x+tu_\alpha )|dt\right|
\emd
Thus to prove the theorem it is enough to establish the inequality
\md1
\big\|\sup_{r>0}|\Gamma_{\alpha ,r,h} f(x)|\big\|_{L^2}\le
C\sqrt{\mu}\|f\|_{L^2}
\emd
By scaling invariance we need only to prove it in the case $h=1$.

\end{section}


\begin{section}{Proof of Theorem}\everypar{} \parskip=11pt

\begin{lemma}\label{L3}
Let $\alpha,\beta\in (0,1)$ be any numbers and $r>0,\,h>0$. The
operator $\Gamma_{\alpha,r,h}f(x)$ defined in \e{Gam} satisfies
estimate
\md1\label{lem1}
|\Gamma_{\alpha ,r,h}f(x)| \le C\left(h
r|\alpha-\beta|+1\right)M_\beta M_{\pi /2}^1f(x), \quad x\in
\ZR^2.
\emd
\end{lemma}

\begin{proof} From \e{int} we have
\md0
 V_r(x_2-x_1\alpha)\le
\zeta_r (x_2-x_1\alpha)
\emd
 Denote
 \md0
 \lambda(x)=2rx|\alpha-\beta|+2
\emd
 and assume
\md1\label{o1}
x_2-x_1\alpha\in (-2^k,2^k),\quad k>\log 1/r .
\emd
for some $k$. Then we shall have $x_2-x_1\alpha\le 2^k$ and so
\md5 \nonumber
\left|\frac{x_2-x_1\beta}{\lambda(x_1)}\right|=
\left|\frac{x_2-x_1\alpha+x_1(\alpha-\beta)}{\lambda(x_1)}\right|\\
\label{Rbig?} \le\left|\frac{x_2-x_1\alpha}{2}\right|+
\frac{1}{2r}\le 2^k,
\emd
which means
\md1\label{o2}
\frac{x_2-x_1\beta}{\lambda(x_1)}\in (-2^k,2^k).
\emd
Hence we conclude that \e{o1} implies \e{o2}. Therefore
\md0
\ZI_{(-2^k,2^k)}(x_2-x_1\alpha)\le
\ZI_{(-2^k,2^k)}\left(\frac{x_2-x_1\beta}{\lambda(x_1)}\right),
\emd
then by definition \e{int} we get
\md5\label{psil}
 V_r(x_2-x_1\alpha)\le
\zeta_r\left(\frac{x_2-x_1\beta}{\lambda(x_1)}\right).
\emd
Now we need the estimate
\md1\label{phil}
\phi_h (x)\le
 \frac{2\left(h r|\alpha-\beta|+1\right)\xi_h(x)}{\lambda(x)}.
\emd
In the case $x>h$ using the estimate $\phi (t)\le \xi (t)/t$ we
get
\md6
\lambda(x)\phi_h (x)=\frac{\lambda(x)}{h}\phi
\left(\frac{x}{h}\right)\le\frac{\lambda(x)}{x}\xi\left(\frac{x}
 {h}\right)=2\left( r|\alpha-\beta|+\frac{1}{x}\right)\xi (\frac{x}{h})\le
\\
 2\left( r|\alpha-\beta|+\frac{1}{h}\right)\xi\left(\frac{x}
 {h}\right)
 =2\left(h r|\alpha-\beta|+1\right)\xi_h(x).
\emd
As $x\le h$ we use $\phi (t)\le \xi (t)$ and so
\md0
\lambda(x)\phi_h (x)\le 2\left(x
r|\alpha-\beta|+1\right)\xi_h(x)\le 2\left(h
r|\alpha-\beta|+1\right)\xi_h(x).
\emd
From\e{phil} and \e{psil} we obtain
\md2
 V_r(x_2-x_1\alpha )\phi_h (x_1){}\le{}
 2\left(h r|\alpha-\beta|+1\right)\xi_h(x_1)
\frac{1}{\lambda(x_1)}\zeta_r\left(\frac{x_2-x_1\beta}{\lambda(x_1)}\right).
\emd
 Finally, using \e{Gam} and \e{int} we conclude
\md6
|\Gamma_{\alpha ,r,h}
f(x)|\le\int_\ZR\xi_h(t_1)dt_1\int_\ZR\frac{1}{\lambda(t_1)}\zeta_r
\left(\frac{t_2-t_1\beta}{\lambda(t_1)}\right)f(x_1+t_1,x_2+t_2)dt_2\\
\le\int_\ZR\xi_h(t_1)M_{\pi /2}^1f(x_1+t_1,x_2+t_1\beta )dt_1\le
M_\beta M_{\pi /2}^1f(x).
\emd

  \end{proof}

For any interval $J=(a,b)$ we denote by $S(J)$ the sector
$\{(x_1,x_2):x_1>0,\, a\le x_2/x_1\le b\}$. For any sector $S$
with angle $\theta $ define by $\gamma S$ ($\gamma
>0$) the sector which has same bisectrix  with $S$ and angle equal
to $\gamma\theta $. Denote by $T_Sf(x)$ the multiplier operator
defined $\widehat {T_Sf}=\ZI_S\widehat f$. For any $c\in (a,b)$ we
consider the restricted strips
\md1\label{str}
S_c(J)=\{(x_1,x_2):x_1>\frac{1}{b-a},\, cx_1-5\le x_2\le cx_1+5\}
\emd

\begin{lemma}\label{L4}
For any complete $\mu $-lacunary set $\Omega $  we have
\md5
\sum_{J\in \seg (\Omega ),\rank (J)\le
\mu-1}\ZI_{S_{p_J}(J)}(x)\le 40,\label{in1}
\\
 \sum_{J=(a,b)\in \seg (\Omega ),\rank (J)= \mu
}\ZI_{S_a(J)}(x)+\ZI_{S_b(J)}(x)\le 12.\label{in2}
\emd
where $p_J\in J$ is the pole of $\Omega $ in $J$.
\end{lemma}
\begin{proof} We shall prove the first inequality. The second one is easier
and can be proved similarly. We consider all intervals $J=(a,b)\in
\seg (\Omega )$ for which $x=(x_1,x_2)\in S_{p_J}(J)$. From the
definition \e{str} we get
\md0
|J|=b-a>\frac{1}{x_1}, \quad \frac{x_2}{x_1}-\frac{5}{x_1}\le p_J
\le\frac{x_2}{x_1}+\frac{5}{x_1}
\emd
and therefore
\md1\label{aJ}
p_J\in [\nu  ,\eta ],\quad |J|>\frac{\eta-\nu}{10}
\emd
We need to prove that, the number of intervals $J\in \seg (\Omega
)$ satisfying \e{aJ} is less than $40$. It can be two cases:
\md4
1)\, J\nsubseteq [\nu ,\eta ],\\ 2)\, J\subset [\nu ,\eta ].
\emd
The set of intervals $J\in \seg (\Omega )$ satisfying 1) and
\e{aJ} forms a monotone sequence
 $J_1\supset J_2\supset \cdots \supset J_l$ of different rank. So we have
 $|J_{k+1}|\le |J_k|/2$. From $p_{J_1}\in [\nu  ,\eta ]$ and
 $J_2\cap [\nu ,\eta ]\neq \varnothing$
 it follows that $|J_2|\le 3\dist(J_2,p_{J_1})\le 3(\eta-\nu )$ and so
 $|J_7|\le 3(\eta-\nu )/32<(\eta-\nu )/10$. Thus we get that $J_i$,
 $i=1,\cdots ,6$ are only intervals satisfying 1) and \e{aJ}.
 Now consider the intervals with  2). We say that such an interval is
 maximal, if there is no other interval with  2) and \e{aJ} in it.
 It is clear the maximal intervals are mutually disjoint, so
 by \e{aJ} their number don't exceed $10$. Each maximal interval can
 be involved at most $3$ higher rank intervals with 2). Thus me
 obtain that the number of all intervals satisfying 2) don't exceed
 $30$. Finally we have got the number of intervals with condition
 \e{aJ} is less that $36$.

\end{proof}


\begin{lemma}\label{L5}
Let $J_1\supset J_2\supset\cdots \supset J_n$ be some sequence of
intervals $J_k=[\alpha_k,\beta_k]\subset (0,1)$ and $p_k\in J_k,\,
k=1,2,\cdots ,n-1 $ satisfies
\md1\label{J}
\frac{|J_{k+1}|}{2}\le \dist(p_k,J_{k+1})\le |J_{k+1}|, \quad
1\le{}k\le{}n-1
\emd
Then for any $\theta \in J_n$ and any function $f\in L^2(\ZR^2)$
we have
\md3 \begin{split} \label{lem2}
\Gamma_{\theta ,R} f\lesssim{}& M_0^2f+M_{\theta }M_{\pi
/2}^1(T_{S_\theta (J_n)}f)\\ {}&\quad{}+\sum_{k=1}^{n-1} M_{p_k}
M_{\pi /2}^1(T_{S_{p_k}(J_k)}f)
\end{split}.
\emd
\end{lemma}

\begin{proof}

  Regard $\theta\in \bigcap J_k$ as fixed. For any $R$
\md1\label{GamR}
\widehat\Gamma_{\theta,R} f(\xi)= \widehat
V_R(\xi_1)\widehat\phi(\xi_2+\xi_1\theta )\widehat f(x).
\emd
Denote
\md1\label{rk}
r_0=0,\quad r_k=\frac{1}{|J_k|}\quad 1\le{}k\le{}m.
\emd
where
\md1\label{mmax}
m=\max\{ k:1\le k\le n,\, 2r_k<R\}.
\emd
To be short in formulas, in some places we shall use notations
$r_{m+1}=R$ and $p_n=\theta $. Defining
\md1\label{Gamk}
\Gamma_kf(x)=\Gamma_{\theta ,r_{k+1}}f(x)-\Gamma_{\theta ,r_k}
f(x)\quad 0\le k\le m,
\emd
 we get
\md1\label{gsum}
\Gamma_{\theta,R} f=\sum_{k=0}^m\Gamma_k f.
\emd
Then by \e{Gam} we have
\md1\label{Gamk1}
\widehat\Gamma_k f(x) = (\widehat V_{r_{k+1}}(\xi_1)- \widehat
V_{r_k}(\xi_1))\widehat\phi(\xi_2+\xi_1\theta ) \widehat f(x),
\quad 0\le k\le m.
\emd
Let us show
\md1\label{supp1}
\supp(\widehat V_{r_{k+1}}(\xi_1)- \widehat
V_{r_k}(\xi_1))\widehat\phi(\xi_2+\xi_1\theta)\subset
S_{p_k}(J_k),\quad 1\le k\le m,.
\emd
Indeed, from \e{psi} and \e{lin} it follows that
\md7
\supp(\widehat V_{r_{k+1}}(\xi_1)- \widehat
V_{r_k}(\xi_1))\widehat\phi(\xi_2+\xi_1\theta)\\ =\{(\xi_1,\xi_2):
r_k\le \xi_1\le 2r_{k+1},\,\, |\xi_2+\xi_1\theta
|<1\}.\label{Vsupp}
\emd
Now let us take any point $(\xi_1,\xi_2)$ from the set \e{Vsupp}.
By \e{rk}
\md1\label{xi1}
\xi_1\ge r_k>\frac{1}{b_k-a_k},\quad
\emd
We have also
\md1
\xi_1\theta -1\le \xi_2\le \xi_1\theta +1
\emd
From \e{J} we easily can get
\md3
b_{k+1}\le p_k+2|J_{k+1}|
\\
a_{k+1}\ge p_k-2|J_{k+1}|
\emd
Therefore using \e{rk} and bound $\xi_1\le 2r_{k+1}$ from
\e{Vsupp} we conclude
\md2
\xi_1\theta +1\le \xi_1b_{k+1} +1\le
\xi_1p_k+2\xi_1|J_{k+1}|+1\le\xi_1p_k+5
\\
 \xi_1\theta -1\ge
\xi_1a_{k+1} -1\ge \xi_1p_k-2\xi_1|J_{k+1}|-1\ge\xi_1p_k-5
\emd
 This inequalities and \e{xi1} implies \e{supp1}.
 From \e{supp1}, \e{Gamk} and \e{Gamk1} it follows that
\md0
\Gamma_k f=\Gamma_k \big(T_{S_{p_k}(J_k)}f\big), \quad
1\le{}k\le{}m
\emd
Hence, using \lem{L3} we conclude
\md1\label{GT}
|\Gamma_k f| {}\lesssim{} (r_{k+1}|\theta-p_k|+1) M_{p_k} M_{\pi
/2}^1\big(T_{S_{p_k}(J_k)}f\big),\quad 1\le{}k\le m.
\emd
Notice also
\md5\label{MM0}
|\Gamma_0 f|\lesssim  M_0^2f.
\emd
If $k=n$, then $|\theta-p_k|=0$. If $k<n$, the by $\theta\in
J_{k+1}\subset J_k$ and \e{J} we have
\md0
|\theta-p_k|\le 2|J_{k+1}|
\emd
The last with \e{rk} implies
\md0
r_{k+1}|\theta-p_k|\le 2
\emd
Hence by \e{GT} we observe
\md0
|\Gamma_k f| {}\lesssim{} M_{p_k} M_{\pi
/2}^1\big(T_{S_{p_k}(J_k)}f\big),\quad 1\le{}k\le m.
\emd
Finally taking account also \e{MM0} we get \lem{L5}.
\end{proof}


\begin{proof}[Proof of Theorem~\ref{Tlac}]

We fix the sets
$\Omega_1\subset\Omega_2\subset\cdots\subset\Omega_{\mu-1}\subset\Omega_\mu=
\Omega $ from definition of N-lacunarity. Fix any angle $\theta
\in \Omega $ and $R>0$. Suppose
\md1
\theta\in \Omega_m\setminus \Omega_{m-1},\hbox { for some } m\le
\mu.
\emd
Denote by $G_k$ the set of all intervals whose vertexes are
neighbor points in $\Omega_k$. We can choose a sequence of
intervals $J_k=[\alpha_k,\beta_k]\in G_k$, $k=1,2,\cdots ,m $ with
poles $p_k\in J_k$ such that
\md0
\theta\in \bigcap_{1\le k\le m} J_k, \qquad \theta=\alpha_m \quad
\text{(or  $\theta=\beta_m$)}
\emd
It is clear that sequence $J_k$ satisfies conditions of \lem{L2}.
Hence,
\md0
\Gamma_{\theta ,R} f\lesssim{} M_0^2f+M_{\theta }M_{\pi
/2}^1(T_{S_\theta(J_m)}f)\\ {}\quad{}+\sum_{k=1}^{m-1} M_{p_k}
M_{\pi /2}^1(T_{S_{p_k}(J_k)}f)
\emd
and therefore,
\md6
|\Gamma_{\theta ,R} f|^2\lesssim{} \mu (|M_0^2f|^2+|M_{\theta
}M_{\pi /2}^1(T_{S_{\theta}(J_m)}f)|^2 +\sum_{k=1}^{m-1 } |M_{p_k}
M_{\pi /2}^1(T_{S_{p_k}(J_k)}f)|^2
\\
\le \mu (|M_0^2f|^2+\sum_{J\in \seg (\Omega )} |M_p M_{\pi
/2}^1(T_{S_{p}(J)}f)|^2
\emd
Hence, using $(2,2)$ bound of strong maximal operators $M_0$ and
$M_p M_{\pi /2}^1$ then \lem{Llac} we obtain
\md7  \label{log}
\int_{\ZR^2}\sup_{\theta\in \Omega , R>0}|\Gamma_{\theta ,R}
f|^2\lesssim \mu (\int_{\ZR^2}|M_0^2f|^2+\sum_{J\in \seg (\Omega
)} \int_{\ZR^2}|M_p M_{\pi /2}^1(T_{S_{p}(J)}f)|^2
\\
\lesssim \mu (\|f\|_{L^2}^2+\sum_{J\in \seg (\Omega )}
\|T_{S_{p}(J)}f)\|_{L^2}^2=\mu (\|f\|_{L^2}^2+\sum_{J\in \seg
(\Omega )} \|\ZI_{S_{p}(J)}\widehat f)\|_{L^2}^2
\\
=\mu (\|f\|_{L^2}^2+\int_{\ZR^2}|\widehat f|^2\sum_{J\in \seg
(\Omega )} \ZI_{S_{p}(J)}\le 40\mu
(\|f\|_{L^2}^2+\int_{\ZR^2}|\widehat f|^2)\lesssim
\mu\|f\|_{L^2}^2.
\emd
Theorem is proved.
\end{proof}

\end{section}


\begin{thebibliography}{99}








\bibitem{Alf}
A.Alfonseca, Strong type inequalities and an almost-orthogonality principle
for families of maximal operators along directions in $\ZR^2 $,
J.London Math.Soc., 67, 2003, No 1, 208-218.
\bibitem{ASV1}
A.Alfonseca, F.Soria, A.Vargas, A remark on maximal operators along directions in
$\ZR^2 $,Math. Res. Lett., 10, 2003,No 1, 41-49.
\bibitem{ASV2}
A.Alfonseca, F.Soria, A.Vargas, An almost-orthogonality principle in $L^2$
for directional maximal functions, Contemp. Math, 2003.
\bibitem{Barr1}
J. Barrionuevo, Estimates for some Kakeya-type maximal operators,
Trans. Amer. Math. Soc.,335,1993, No 2, 667-682
\bibitem{Barr}
J.Barrionuevo, A note on the Kakeya maximal operator, Math. Res.
Letters,3, 1996, No1, 61-65.
\bibitem{Car}
A. Carbery, Differentiation in lacunary directions and an extension of the
            Marcinkiewicz multiplier theorem,Ann. Inst. Fourier (Grenoble),
            38, 1988, No1,157-168.
\bibitem{Cor}
A. C\'{o}rdoba, The Kakeya maximal function and the spherical summation
            multipliers, Amer. J. Math., 99,1977, No 1, 1-22.
\bibitem{Chr}
M.Christ, Examples of singular maximal functions unbounded on $L\sp p$,
Conference on Mathematical Analysis (El Escorial, 1989),Publ. Mat.,
35, 1991, No 1, 269-279.
\bibitem{CoFe}
A.C\'{o}rdaba and R.Fefferman, On differentiation of integrals, Proc. Nat. Acad.
of Sci USA, 74, 1977, No 2, 423-425.
\bibitem{CoFe1}
   A.C\'{o}rdaba and R.Fefferman, On the equivalence between the boundedness of certain classes of
            maximal and multiplier operators in Fourier analysis,
Proc. Nat. Acad. Sci. U.S.A., 74, 1977, No 2, 423-425.
\bibitem{DV}
J.Duoandikoetxea, A.Vargas, Directional operators and radial functions on
the plane, Ark. Mat., 33, 1995, No 2, 281-291.
\bibitem{Har}
K.Hare, Maximal operators and Cantor sets, Canad. Math. Bull., 43, 2000,
330-342.
\bibitem{HR}
K.Hare, F.Ricci, Maximal functions with polynomial densities in lacunary
directions, Trans. Amer. Math. Soc., 335, 2003, No 3, 1135-1144.
\bibitem{HRo1}
K.Hare, J.-O. R{\"o}nning
Applications of generalized Perron trees to maximal functions
and density bases, J. Fourier Anal. Appl., 4, 1998, No 2, 215-227.
\bibitem{HRo2}
K.Hare, J.-O. R{\"o}nning, The size of ${\rm Max}(p)$ sets and density bases,
J. Fourier Anal. Appl., 8, 2002, No 3, 259-268.
\bibitem{KaLa}
G.A.Karagulyan, M.Lacey, An estimate of the maximal operators
associated with generalized lacunary sets, Contemporary Math.
Anal., 2005, No 1.

\bibitem{Katz}
N.H.Katz, A counterexample for maximal operators over a Cantor set of
directions, Math. Res. Lett., 3, 1996, No 4, 527-536.
\bibitem{Katz1}
N.H.Katz, Remarks on maximal operators over arbitrary sets of directions,
Bull. London Math. Soc., 31, 1999, No 6, 700-710.
\bibitem{Katz2}
N.H.Katz, Maximal operators over arbitrary sets of directions,
Duke Math. J., 97, 1999, No 3, 67-79.
\bibitem{NSW}
A.Nagel, E.M.Stein, S.Wainger, Differentiation in lacunary directions.
Duke Math. J., 97, 1979, No. 1, 67-79.
\bibitem{SS}
P.Sjogren and P.Sjolin, Littlewood-Paley decompositions and
Fourier multipliers with singularities on certain sets, Ann. Inst. Fourier
(Grenoble),31, 1981, No 1, 157-175.
\bibitem{Str1}
J.-O.Str{\"o}mberg, Weak estimates for maximal functions with
rectangles in certain directions, Arkiv. for Mat,15, 1978,
229-240.
\bibitem{Str2}
J.-O.Str{\"o}mberg, Maximal functions associated to rectangles with uniformaly
distributed directions, Ann. of Math.,107, 1976, 399-402.
\bibitem{Var}
A.Vargas, A remark on a maximal function over a Cantor set of
directions, Rend. Circ. Mat. Palermo (2), 44, 1995, No 2, 273-282.
\bibitem{Wai}
S.Wainger, Applications of Fourier transforms to averages over
lower-dimensional sets, Harmonic analysis in Euclidean spaces
(Proc. Sympos. Pure Math.,Williams Coll., Williamstown, Mass., 1978), Part 1,
Proc. Sympos. Pure Math., XXXV, Part, Amer. Math. Soc., Providence, R.I.,
1979, 85-94
\end{thebibliography}
\end{document}